# Some remarks about the minimal polynomials of $4\sin^2\left(\frac{\pi}{n}\right)$.

Johann Cigler

johann.cigler@univie.ac.at

**Abstract**

We show that the values of the minimal polynomials $\phi_n(x)$ of $4\sin^2\left(\frac{\pi}{n}\right)$ for $x \in \{0,1,2,3,4\}$ are intimately related to the prime factorization of $n$.

**1. Introduction**

Let $L_n(x)$ denote the Lucas polynomials defined by $L_n(x) = xL_{n-1}(x) - L_{n-2}(x)$ with initial values $L_0(x) = 2$ and $L_1(x) = x$. In [1] the polynomials $Z_n(x) = 2 - L_n(2-x)$ were studied which are closely related to the spread polynomials of Norman J. Wildberger and are characterized by $Z_n(4\sin^2\theta) = 4\sin^2(n\theta)$. They satisfy $Z_n(x) = \prod_{d|n} \Phi_d(x)$ with polynomials $\Phi_d(x) \in \mathbb{Z}[x]$ of degree $\deg \Phi_d = \varphi(d)$, where $\varphi(n)$ is Euler's totient function. The polynomials $\Phi_n(x)$ are given by $\Phi_1(x) = x = \phi_1(x)$, $\Phi_2(x) = 4 - x = \phi_2(x)$ and $\Phi_n(x) = \phi_n^2(x)$ for $n \geq 3$, where $\phi_n(x)$ is the minimal polynomial of $4\sin^2\left(\frac{\pi}{n}\right)$, normalized such that the constant term is positive.

For example,

$$\left(L_n(x)\right)_{n\geq 0} = \left(2, x, -2+x^2, -3x+x^3, 2-4x^2+x^4, 5x-5x^3+x^5, -2+9x^2-6x^4+x^6, \cdots\right),$$

$$\left(Z_n(x)\right)_{n\geq 1} = \left(x, 4x-x^2, 9x-6x^2+x^3, 16x-20x^2+8x^3-x^4, \cdots\right),$$

$$\left(\phi_n(x)\right)_{n\geq 1} = \left(x, 4-x, 3-x, 2-x, 5-5x+x^2, 1-x, 7-14x+7x^2-x^3, 2-4x+x^2, 3-9x+6x^2-x^3, \cdots\right).$$

For $x \in \{0,1,2,3,4\}$ the sequence $\left(Z_n(x)\right)_{n\geq 0}$ is periodic:

$$\left(Z_n(0)\right)_{n\geq 0} = (0,0,0,\cdots),$$

$$\left(Z_n(1)\right)_{n\geq 0} = (0,1,3,4,3,1,0,1,3,4,3,1,0,1,3,4,3,1,\cdots),$$

$$\left(Z_n(2)\right)_{n\geq 0} = (0,2,4,2,0,2,4,2,0,2,4,2,0,2,4,2,\cdots),$$

$$\left(Z_n(3)\right)_{n\geq 0} = (0,3,3,0,3,3,0,3,3,0,3,3,0,3,3,\cdots),$$

$$\left(Z_n(4)\right)_{n\geq 0} = (0,4,0,4,0,4,0,4,0,4,0,4,0,4,0,4,\cdots).$$

It turns out that the corresponding values of $\phi_n(x)$ for $n > 2$ are signed primes as the following examples show.

$$(\phi_n(0))_{n\geq 1} = (0,4,3,2,5,1,7,2,3,1,11,1,13,1,1,2,17,1,19,1,1,1,23,\cdots),$$

$$(\phi_{2n}(4))_{n\geq 1} = (0,-2,-3,2,5,1,-7,2,-3,1,-11,1,13,1,1,2,17,1,-19,1,1,1,-23,\cdots),$$

$$(\phi_{3n}(3))_{n\geq 1} = (0,-2,3,-2,-5,1,7,-2,3,1,-11,1,13,1,1,-2,-17,1,19,1,1,1,-23,\cdots),$$

$$(\phi_{4n}(2))_{n\geq 1} = (0,-2,-3,2,5,1,-7,2,-3,1,-11,1,13,1,1,2,17,1,-19,1,1,1,-23,\cdots),$$

$$(\phi_{6n}(1))_{n\geq 1} = (0,-2,-3,-2,-5,1,7,-2,-3,1,-11,1,13,1,1,-2,-17,1,19,1,1,1,-23,\cdots).$$

The purpose of this note is a clarification of this situation.

## 2. A formula for the minimal polynomials

In [1], (5.3) it has been shown that the minimal polynomials with positive constant term are given by

(1) $$\phi_n(x) = \frac{C_n(\lambda(x))}{(\lambda(x))^{\frac{\varphi(n)}{2}}}$$

for $n \geq 3$, where $C_n(x) = \prod_{\substack{1 \leq k \leq n \\ k \perp n}} (x - \zeta_n^k)$ with $\zeta_n = e^{\frac{2\pi i}{n}}$ is the $n$-th cyclotomic polynomial, $\lambda(x) = \frac{2 - x + \sqrt{x^2 - 4x}}{2}$, and $\varphi(n)$ is Euler's totient function.

In order to make the paper self-contained let me recall the proof of (1).

For the Lucas polynomials $L_n(x)$ Binet's formula gives $L_n(x) = \alpha(x)^n + \beta(x)^n$ where $\alpha(x) = \frac{x + \sqrt{x^2 - 4}}{2}$ and $\beta(x) = \frac{x - \sqrt{x^2 - 4}}{2} = \frac{1}{\alpha(x)}$ are the roots of $z^2 = xz - 1$.

For $n \geq 3$ the cyclotomic polynomial $C_n(x)$ is a monic symmetric polynomial with integer coefficients of even degree $\varphi(n)$. Therefore $\frac{C_n(x)}{x^{\frac{\varphi(n)}{2}}}$ is a sum of terms $x^k + \frac{1}{x^k}$ with integer coefficients. Since $\alpha(x)^k + \frac{1}{\alpha(x)^k} = L_k(x)$ we get $\frac{C_n(\alpha(x))}{\alpha(x)^{\frac{\varphi(n)}{2}}} = \sum_{k=0}^{\frac{\varphi(n)}{2}} c_k L_k(x) \in \mathbb{Z}[x]$ with integer coefficients $c_k$. Since $\alpha(2\cos\theta) = \frac{2\cos\theta + \sqrt{(2\cos\theta)^2 - 4}}{2} = \cos\theta + i\sin\theta = e^{i\theta}$ we get $C_n\left(\alpha\left(2\cos\left(\frac{2\pi}{n}\right)\right)\right) = 0$. Thus $\frac{C_n(\alpha(x))}{\alpha(x)^{\frac{\varphi(n)}{2}}} \in \mathbb{Z}[x]$ is a monic polynomial of degree $\frac{\varphi(n)}{2}$ with root $2\cos\left(\frac{2\pi}{n}\right)$.

Therefore it coincides with the minimal polynomial

(2) $$\psi_n(x) = \prod_{0<j<\frac{n}{2}, j\perp n} \left(x - 2\cos\left(\frac{2\pi j}{n}\right)\right)$$

of $2\cos\left(\frac{2\pi}{n}\right)$. Observing that $2 - 2\cos\left(\frac{2\pi k}{n}\right) = 4\sin^2\left(\frac{k\pi}{n}\right)$ we get by changing $x \to 2-x$

$$\psi_n(2-x) = \prod_{0<j<\frac{n}{2}, j\perp n} \left(2 - x - 2\cos\left(\frac{2\pi j}{n}\right)\right) = \prod_{0<j<\frac{n}{2}, j\perp n} \left(4\sin^2\left(\frac{j\pi}{n}\right) - x\right).$$

Thus $\phi_n(x) = \prod_{0<j<\frac{n}{2}, j\perp n} \left(4\sin^2\left(\frac{j\pi}{n}\right) - x\right) = \psi_n(2-x) = \dfrac{C_n(\lambda(x))}{(\lambda(x))^{\frac{\varphi(n)}{2}}}$

with $\lambda(x) = \alpha(2-x)$.

For example, $C_3(x) = 1 + x + x^2$, $\dfrac{C_3(x)}{x} = \left(x + \dfrac{1}{x}\right) + 1$, $\psi_3(x) = x+1$, $\phi_3(x) = 2 - x + 1 = 3 - x$,

$C_4(x) = 1 + x^2$, $\dfrac{C_4(x)}{x} = x + \dfrac{1}{x}$, $\psi_4(x) = x$, $\phi_4(x) = 2 - x$,

$C_5(x) = 1 + x + x^2 + x^3 + x^4$, $\dfrac{C_5(x)}{x^2} = 1 + \left(x + \dfrac{1}{x}\right) + \left(x^2 + \dfrac{1}{x^2}\right)$,

$\psi_5(x) = 1 + x + (x^2 - 2) = x^2 + x - 1$, $\phi_5(x) = (2-x)^2 + (2-x) - 1 = x^2 - 5x + 5$.

### 3. Some background about cyclotomic polynomials

We need some well-known results about cyclotomic polynomials:

**Lemma 1**

*For a prime number $p$*

(3) $$C_{pn}(x) = C_n(x^p)$$

*if $p | n$ and*

(4) $$C_{pn}(x) = \dfrac{C_n(x^p)}{C_n(x)}$$

*if $p$ and $n$ are relatively prime, which will be noted as $p \perp n$.*

*For odd numbers $n$ we get*

(5) $$C_{2n}(x) = C_n(-x).$$

**Proof**

Note that

(6)
$$\varphi(pn) = p\varphi(n) \quad \text{if } p \mid n,$$
$$= (p-1)\varphi(n) \quad \text{if } p \perp n$$

For $p \mid n$ $C_{pn}(x)$ is irreducible with degree $\varphi(pn) = p\varphi(n)$ with root $\zeta_{pn}$ and $C_n(x^p)$ has degree $p\phi(n)$ and root $\zeta_{np}$. Therefore it must coincide with $C_{pn}(x)$.

For $p \perp n$ $C_{pn}(x)$ is irreducible with degree $\phi(pn) = (p-1)\varphi(n)$ with root $\zeta_{pn}$.

$\dfrac{C_n(x^p)}{C_n(x)}$ has the same root and the same degree.

For odd $n$ $C_{2n}(x)$ and $C_n(-x)$ have both degree $\varphi(n)$ and $C_{2n}(x)$ is irreducible with root $\zeta_{2n}$. Since $C_n(-x)$ has the same root $-\zeta_n = \zeta_{2n}$ the polynomials must coincide.

The following result seems to be well known.

**Lemma 2**

For $n > 1$

(7)
$$C_n(1) = v(n),$$

where $v(p^k) = p$ for powers $k \geq 1$ of a prime $p$ and $v(n) = 1$ else.

**Proof**

For a prime $p$ $C_p(x) = \dfrac{x^p - 1}{x - 1}$ implies $C_p(1) = p = v(p)$.

We use induction with respect to the index $n$.

(7) is true for $n = 2$.

Let (7) be true for $m < n$,

If $n$ is a prime then (7) is true.

If not, let $p$ be a prime factor of $n$.

If $n = p^k$ then by (3) $C_{p^k}(1) = C_{p^{k-1}}(1) = \cdots = C_p(1) = p = v(p^k)$.

If not, then $n = p^k m$ with $m \perp p$ and $m > 1$. By (3) $C_n(1) = C_{pm}(1)$ and by (4)

$$C_{pm}(1) = \dfrac{C_m(1^p)}{C_m(1)} = 1 = v(pm) = v(n).$$

**Lemma 3**

Let $w(n,x) = \dfrac{C_n(x)}{x^{\frac{\varphi(n)}{2}}}$. Then for a prime $p$ we get

(8) $$w(pn, x) = w(n, x^p)$$

if $p \mid n$ and

(9) $$w(pn, x) = \dfrac{w(n, x^p)}{w(n, x)}$$

for $p \perp n$.

**Proof**

By (6) we get for $p \mid n$

$$w(pn, x) = \dfrac{C_{pn}(x)}{x^{\frac{\varphi(pn)}{2}}} = \dfrac{C_n(x^p)}{(x^p)^{\frac{\varphi(n)}{2}}} = w(n, x^p) \text{ and for } p \perp n$$

$$w(pn, x) = \dfrac{C_{pn}(x)}{x^{\frac{\varphi(pn)}{2}}} = \dfrac{C_n(x^p)}{C_n(x)} \cdot \dfrac{1}{x^{\frac{(p-1)\varphi(n)}{2}}} = \dfrac{w(n, x^p)}{w(n, x)}.$$

**3. Results**

**Theorem 1**

For $n > 2$

(10) $$\phi_n(0) = v(n).$$

**Proof**

This follows from (1) since $\lambda(0) = 1$.

**Remark**

It should be noted that Theorem 1 also follows from the formula $\displaystyle\prod_{\substack{0<k<n \\ k \perp n}} 2\sin\left(\dfrac{k\pi}{n}\right) = v(n)$ which has been proved in [2],(4).

**Theorem 2**

*For $n \geq 3$*

$\phi_{2n}(4) = v(n)$ if $v(n) \equiv 1 \bmod 4$ or $v(n) = 2$ and $\phi_{2n}(4) = -v(n)$ if $v(n) \equiv -1 \bmod 4$.

**Proof**

We know that $\phi_n(4) = w(n,-1) = (-1)^{\frac{\varphi(n)}{2}} C_n(-1)$.

We first show that $w(2n+1,-1) = \varepsilon_{2n+1}$ with $\varepsilon_n = -1$ if $v(n) \equiv -1 \pmod 4$ and $\varepsilon_n = 1$ else.

For an odd prime $p$ we have $C_p(-1) = 1$. For $p \equiv 1 \pmod 4$ we have $\varphi(p) \equiv 0 \pmod 4$ and for $p \equiv -1 \pmod 4$ we get $\varphi(p) \equiv 2 \pmod 4$. For $p = 2$ we get $w(2,-1) = 0$.

Thus for primes $p > 2$ we have $w(p,-1) = \varepsilon_p$. By Lemma 2 then also $w(p^k,-1) = \varepsilon_p$. If $2n+1 = p^k m$ with $p \perp m$ then $w(2n+1,-1) = w(pm,-1) = \dfrac{w(m,(-1)^p)}{w(m,-1)} = 1$.

Thus for odd $n > 2$ $w(n,-1) = -1$ if $\varepsilon_n = -1$ and $w(n,-1) = 1$ else.

If $n = 2m$ we have either $n = 2^k$ with $w(2^k,-1) = w(2^{k-1},1) = 2 = v(2^k)$ or $n = 2^k m$ with odd $m > 1$. Then $w(n,-1) = w(2^k m,-1) = w(2m,-1) = \dfrac{w(m,1)}{w(m,-1)} = \varepsilon_m = \varepsilon_m v(m)$. This gives Theorem 2.

Because $\phi_{4n}(2) = w(4n, i) = w(2n,-1) = \phi_{2n}(4)$ we get

**Theorem 3**

*For $n \geq 3$*

$\phi_{4n}(2) = v(n)$ *if* $v(n) \equiv 1 \bmod 4$ *or* $v(n) = 2$ *and* $\phi_{4n}(2) = -v(n)$ *if* $v(n) \equiv -1 \bmod 4$.

**Theorem 4**

*For $n \geq 3$*

$\phi_{3n}(3) = v(n)$ *if* $v(n) \equiv 1 \bmod 3$ *or* $v(n) = 3$ *and* $\phi_{3n}(3) = -v(n)$ *if* $v(n) \equiv -1 \bmod 3$.

**Proof**

$\lambda(3) = \dfrac{-1+\sqrt{-3}}{2} = \zeta_3 = \omega$ and $\varphi_n(3) = \dfrac{C_n(\omega)}{\omega^{\frac{\varphi(n)}{2}}} = w(n,\omega)$.

Let $\delta_n = -1$ if $v(n) \equiv -1 \pmod 3$ and $\delta_n = 1$ else.

We first show that $w(n,\omega) = \delta_n$ for $n > 3$ with $n \perp 3$.

Since $C_4(x) = x^2 + 1$ we get $w(4,\omega) = \dfrac{C_4(\omega)}{\omega^{\frac{\varphi(4)}{2}}} = \dfrac{\omega^2 + 1}{\omega} = -\dfrac{\omega}{\omega} = -1 = \delta_4$.

For a prime $p = 6n+1$ we get $w(p,\omega) = \dfrac{\omega^p - 1}{\omega - 1} \dfrac{1}{\omega^{\frac{\varphi(p)}{2}}} = 1 = \delta_p$, for $p = 6n+5$ we get

$$w(p,\omega) = \dfrac{\omega^p - 1}{\omega - 1} \dfrac{1}{\omega^{\frac{\varphi(p)}{2}}} = \dfrac{\omega^2 - 1}{\omega - 1} \dfrac{1}{\omega^2} = -1 = \delta_p.$$

Now we use induction. Suppose that $w(m,\omega) = \delta_m$ for $3 < m < n$ with $m \perp 3$. If $n$ is a prime then $w(n,\omega) = \delta_n$. If not, let $n = pm$ for a prime $p$.

If $n = p^k$ then $w\left(pp^{k-1}, \omega\right) = w\left(p^{k-1}, \omega^p\right)$.

For $p \equiv 1 \pmod 3$ we get $w\left(p^k, \omega\right) = w(p,\omega) = \delta_p = 1$.

For $p \equiv 2 \pmod 3$ we get $w\left(p^k, \omega\right) = w\left(p^{k-1}, \omega^2\right) = w\left(p, \omega^{2k}\right) = w(p,\omega) = -1 = \delta_{p^k}$, because $\omega^2 = \bar\omega$ and therefore $w(x,\omega^2)$ is the complex conjugate of $w(x,\omega)$.

Let now $n = p^k m$ with $p \perp m$. Then we get

$$w(n,\omega) = w\left(pm, \omega^{k-1}\right) = w(pm,\omega) = \dfrac{w(m,\omega)}{w(m,\omega)} = 1 = \delta_n.$$

Finally we get for $n = p^k$ with $p \perp 3$ that $w(3n,\omega) = \dfrac{w(n,1)}{w(n,\omega)} = \delta_n v(n)$.

For $n = 3^k$ we get $w(3^k,\omega) = w(3^{k-1},1) = v(3^k) = 3$.

If $n = 3^k m$ with $m \perp 3$ then $w(3^k m, \omega) = w(3m,1) = 1 = v(n)$.

**Theorem 5**

*For $n \geq 3$*

$\phi_{6n}(1) = v(n)$ if $v(n) \equiv 1 \bmod 3$ and $\phi_{6n}(1) = -v(n)$ if $v(n) \equiv -1 \bmod 3$ or if $v(n) = 2$ or $v(n) = 3$.

**Proof**

Note that $\phi_n(1) = w(n,\sigma)$ with $\sigma = \dfrac{1 + \sqrt{-3}}{2} = -\omega^2 = \zeta_6$.

The result follows from

(11) $$w(6n, \sigma) = w(3n, \omega)$$

for $n \geq 2$ with the exception that for $n = 3^k$

(12) $$w(6 \cdot 3^k, \sigma) = -w(3^{k+1}, \omega).$$

Formula (11) holds for even $n$ because

$$w(6n, \sigma) = w(2 \cdot 3n, \sigma) = w(3n, \sigma^2) = w(3n, \omega) \text{ if } 2 \mid n.$$

For $n = 3^k$ we get from (5)

$$w(2 \cdot 3^{k+1}, \sigma) = \frac{C(2 \cdot 3^{k+1}, \sigma)}{\sigma^{\frac{\varphi(2 \cdot 3^{k+1})}{2}}} = \frac{C(3^{k+1}, -\sigma)}{\sigma^{\frac{\varphi(2 \cdot 3^{k+1})}{2}}} = \frac{C(3^{k+1}, \omega^2)}{(-\omega^2)^{3^k}} = -w(3^{k+1}, \omega^2) = -w(3^{k+1}, \omega) = -3.$$

For a prime $p$ we get from (5)

$$w(6p^{n+1}, \sigma) = \frac{C(2 \cdot 3p^{n+1}, \sigma)}{\sigma^{\frac{\varphi(6p^{n+1})}{2}}} = \frac{C(3p^{n+1}, -\sigma)}{\sigma^{p^n(p-1)}} = \frac{C(3p^{n+1}, \omega^2)}{(-\sigma)^{p^n(p-1)}} = \frac{C(3p^{n+1}, \omega^2)}{(\omega^2)^{p^n(p-1)}}$$

$$= w(3p^{n+1}, \omega^2) = w(3p^{n+1}, \omega) = \delta_p p.$$

There remains the case $v(n) = 1$. Since the formula holds for even $n$ we can assume that $n$ is a product of at least 2 different odd primes.

Let us first consider the case $n \equiv 0 \pmod 3$.

$$w(6n, \sigma) = \frac{C_{6n}(\sigma)}{\sigma^{\frac{\varphi(6n)}{2}}} = \frac{C_{3n}(-\sigma)}{\sigma^{\frac{3\varphi(n)}{2}}} = \frac{C_{3n}(\omega^2)}{(-\omega^2)^{\frac{3\varphi(n)}{2}}} = \frac{C_{3n}(\omega^2)}{(-\omega^2)^{\frac{\varphi(3n)}{2}}} = w(3n, \omega^2) = w(3n, \omega).$$

If $n \not\equiv 0 \pmod 3$ we get

$$w(6n, \sigma) = \frac{C_{6n}(\sigma)}{\sigma^{\frac{\varphi(6n)}{2}}} = \frac{C_{3n}(-\sigma)}{\sigma^{\varphi(3n)}} = \frac{C_{3n}(\omega^2)}{(-\omega^2)^{\frac{\varphi(3n)}{2}}} = w(3n, \omega^2) = w(3n, \omega).$$

**Addendum**

After posting the first version of this note I found reference [3], where already some of the results about cyclotomic polynomials have been proved.

**References**

[1] Johann Cigler and Hans-Christian Herbig, Factorization of spread polynomials, arXiv: 2412.18958 (2024)
[2] Peter Luschny and Stefan Wehmeier, The LCM(1,2,…,n) as a product of sine values, sampled over the points in Farey sequences, arXiv:0909.1838 (2009)
[3] Bartolomiej Bzdega, Andres Herrera-Poyatos and Pieter Moree, Cyclotomic polynomials at roots of unity, arXiv:1611.06783 (2017)

9